\documentclass{article}
\usepackage{eurosym}
\usepackage{amssymb}
\usepackage{amsfonts}
\usepackage{amsmath}
\usepackage{graphicx}
\usepackage{hyperref}
\usepackage{subcaption}
\usepackage{tikz}
\usepackage{float}

\newtheorem{theorem}{Theorem}

\newtheorem{condition}[theorem]{Condition}

\newtheorem{lemma}[theorem]{Lemma}

\newtheorem{remark}[theorem]{Remark}

\begin{document}

\author{ Pelin G. Geredeli \thanks{%
email address: peling@iastate.edu.} \\
Department of Mathematics\\
Iowa State University, Ames-IA, USA}
\title{A Time Domain Approach for the Exponential Stability of a
Nondissipative Linearized Compressible Flow-Structure PDE System}
\maketitle

\begin{abstract}
This work is motivated by a longstanding interest in the long time behavior of
flow-structure interaction (FSI) PDE dynamics. Such coupled PDE systems are ubiquitous
in modeling of various natural phenomena, and in particular have many applications in fluid dynamics, aeroelasticity and
gas dynamics. We consider a linearized compressible flow
structure interaction (FSI) PDE model with a view of analyzing the stability properties
of both the compressible flow and plate solution components. In our earlier  work, we gave an answer in the affirmative to question of uniform
stability for finite energy solutions of said compressible flow-structure system, by means of a ``frequency domain'' approach. However, the frequency domain method of proof in  that work is not ``robust'' (insofar as we can see), when one wishes to study longtime behavior of solutions of compressible flow-structure PDE models which track the appearance of the ambient state onto the boundary interface (i.e., $\kappa =1$ in flow PDE component (2)). Nor is a frequency domain approach in this earlier work availing when one wishes to consider the dynamics, in long time, of solutions to
physically relevant nonlinear versions of the compressible flow-structure PDE system under present consideration (e.g., the Navier-Stokes nonlinearity in the PDE flow
component, or a nonlinearity of Berger/Von Karman type in the plate equation). Accordingly, in the present work, we operate in the time domain by way of obtaining the necessary energy estimates which culminate in an alternative proof  for the uniform stability of finite energy compressible flow-structure solutions. This novel time domain proof will be used in our forthcoming paper which addresses the
existence of compact global attractors for the corresponding nonlinear
coupled system in which the material derivative -- which incorporates the aforesaid ambient state -- of the interaction surface will be taken into account. Since there is a need to avoid steady states in
our stability analysis, as a prerequisite result, we also show here that zero is an eigenvalue for the generators of flow-structure systems, whether the material derivative term be absent or present. Moreover, we provide a clean characterization of the (one dimensional) zero eigenspace, with or without material derivative, under an appropriate assumption on the underlying ambient vector field.

 \vskip.3cm
\noindent \textbf{Key terms:} Flow-structure interaction, compressible
flows, uniform stability
\end{abstract}


\section{\protect\bigskip Introduction}

The linearized compressible flow-structure PDE model which we will consider arises
in the context of the design of various engineering systems and the study of
gas dynamics. This coupled system describes the interaction between a plate
and a given compressible gas flow. In contrast to incompressible fluid flows,
wherein the fluid density is assumed to be a constant, compressible gas flow models will contain an additional fluid density variable,
and involve other state spaces. For further details the reader is referred to \cite{Chu2013-comp, agw, material, pelin-george}.

The presence
of the density (pressure) equations in compressible cases will tend to  make the analysis
quite different than that for incompressible flows; in particular, one must deal with the extra density
(pressure) solution variable. Moreover, and intrinsic to the problem under consideration, linearization of
the compressible Navier Stokes equation occurs around a rest state that
contains an arbitrary ambient vector field $\mathbf{U}$. Qualitative properties of this model -- i.e., wellposedness and long term analysis in the sense of global attractors (in the presence of the von Karman plate nonlinearity) were firstly analyzed in \cite{Chu2013-comp}, in the case $\mathbf{U}=0.$ However, the case $\mathbf{U}\neq 0$ was recognized to be challenging, since the presence of a nonzero ambient field introduces problematic terms such as $\mathbf{U} \cdot \nabla p,$ (where $p$ is the pressure variable), a term which is ``unbounded'' with respect to the underlying finite energy space of wellposedness. 

With respect to compressible flow-structure PDE systems with underlying nonzero ambient terms: A positive answer to the wellposedness question was given in \cite{agw} in the case $\mathbf{U}\neq 0$; exponential stability of finite energy solutions to this model was shown in \cite{pelin-george}, again in the case $\mathbf{U}\neq 0$. With a view of handling the aforesaid troublesome term $\mathbf{U} \cdot \nabla p$, a frequency domain approach  is invoked in \cite{pelin-george}. By way of appropriately estimating $\mathbf{U} \cdot \nabla p$ in \cite{pelin-george}, as a static (and not time dependent term), the frequency domain approaches allows for an appropriate decomposition of static Stokes flow, and an eventual invocation of the nonsmooth domain version of  the Agmon-Douglis-Nirenberg Theorem; see p. 75 of \cite{dauge}. In addition to dealing with unbounded term $\mathbf{U} \cdot \nabla p$, a large part of the work in \cite{pelin-george} is devoted to a spectral analysis of the compressible flow-structure generator, by way of ultimately invoking the wellknown resolvent criteria for exponential decay in \cite{huang} and \cite{pruss}.

Although the methodology set forth in \cite{pelin-george} is effective in establishing exponential stability for solutions of the compressible flow-structure PDE model (2)-(4) below, with $\kappa =0$, it's use seems limited when dealing with (2)-(4) when the physically relevant material derivative term is present (i.e., $\kappa =1$; see \cite{material} for the modeling aspects of this problem; also \cite{dowell1}). In particular, since the presence of the material derivative term $\mathbf{U}\cdot \nabla w$ on the boundary interface constitutes an unbounded perturbation of the compressible flow-structure PDE system, a necessary spectral analysis for $\kappa =1$, analogous to that in \cite{pelin-george}, is problematic. In addition, the critical frequency domain estimates which were obtained in \cite{pelin-george} for the linear problem do not lend themselves readily to adaptation so as to handle nonlinear versions of (2)-(4), versions in which the Navier-Stokes or von Karman plate nonlinearities are present. Consequently, the problem of analyzing long time behavior for nonlinear compressible flow-structure PDE systems -- particularly in the sense of global attractors -- must be undertaken in the "time domain" rather than the "frequency domain". 

 Accordingly, the principal contribution of the present manuscript is to give an alternative proof for the exponential stability of the solutions via a certain multiplier method in the  time domain. In particular, our main (gradient type) multiplier is based upon the solution of a certain Neumann problem, a solution which is sufficiently smooth, even considering the unavoidable boundary interface singularities; see \cite{Dauge_1} and \cite{JK}. We should also note that the application of this multiplier is practicable and convenient, due to the characterization (compatibility condition) of the stabilizable finite energy space $[Null(\mathcal{A})]^{\bot }$ where $\mathcal{A}$ is the semigroup generator of the system.  

As we said, besides being of intrinsic interest in its own right, this manuscript will also serve as a blue print for our forthcoming work which will address the existence of compact attractors for nonlinear flow-structure PDE interactions, in which the material derivative appears in the normal component boundary condition of the compressible flow variable (i.e., $\kappa = 1$). This material derivative term on the boundary interface represents an unbounded perturbation of the compressible flow-structure semigroup generator, in addition to the presence of a nonlocal nonlinearity (von Karman or Berger). With this future work in mind, in this manuscript, we also analyze the Null space of the generator of the system, in the presence or absence of the material derivative on the interaction
surface $\Omega$.  Again, our main reason for doing this is to avoid finite energy initial data which gives rise to steady states. In the course of proof, which partly involves a necessary assumption on the ambient field $\mathbf{U}$, we show that, despite the unboundedness presented by the material derivative perturbation, the Null Space of the ``material derivative'' generator $\mathcal{A}_{1}$, $(\kappa=1)$ will coincide with Null Space of the ``material derivative free'' generator $\mathcal{A}_{0}$, $(\kappa=0).$

\bigskip


In what follows we provide the PDE description of the interaction model
under the study. Let the flow domain $\mathcal{O}$ be a bounded subset of $%
\mathbb{R}^{3}$, with boundary $\partial \mathcal{O}$. Moreover, $\partial 
\mathcal{O}=\overline{S}\cup \overline{\Omega }$, with $S\cap \Omega
=\emptyset $, and with (structure) domain $\Omega \subset \mathbb{R}^{3}$
being a \emph{flat} portion of $\partial \mathcal{O}$. In particular, assume that $%
\partial \mathcal{O}$ has the following specific configuration: 
\begin{equation}
\Omega \subset \left\{ x=(x_{1,}x_{2},0)\right\} \,\text{\ and \ surface }%
S\subset \left\{ x=(x_{1,}x_{2},x_{3}):x_{3}\leq 0\right\} \,.  \label{geo}
\end{equation}%
and $\ \mathbf{n}(\mathbf{x})$ denotes the unit outward normal
vector to $\partial \mathcal{O}$ where $\left. \mathbf{n}\right\vert _{\Omega
}=[0,0,1]$. We also impose additional following geometrical conditions on $\mathcal{O}:$

\begin{condition}

\label{cond} The flow domain $\mathcal{O}$ should be a curvilinear polyhedral
domain --i.e., $\mathcal{O}$ has a finite set of smooth \ edges and corners;
see, \cite{Dauge_3} -- which satisfies the following assumptions: 

\begin{enumerate}
\item   Each corner of the boundary  $\partial \mathcal{O}$ { -if any- } is diffeomorphic to a convex cone, 
\item  Each point on an edge of the boundary  $%
\partial \mathcal{O}$ \text{ is diffeomorphic to a wedge with opening} $<\pi. 
$ 

\end{enumerate}

\end{condition}

\begin{remark}
The point of making such assumptions on the geometry of domain $\mathcal{O}$  is that they will allow for our application of elliptic results for solutions of second order boundary value problems on corner domains; see \cite{Dauge_1}. In particular, these results will be invoked in our time dependent multiplier method by way of proving our uniform stability result.
\end{remark}
\noindent The usual, familiar geometries can be considered as
in Figure 1.

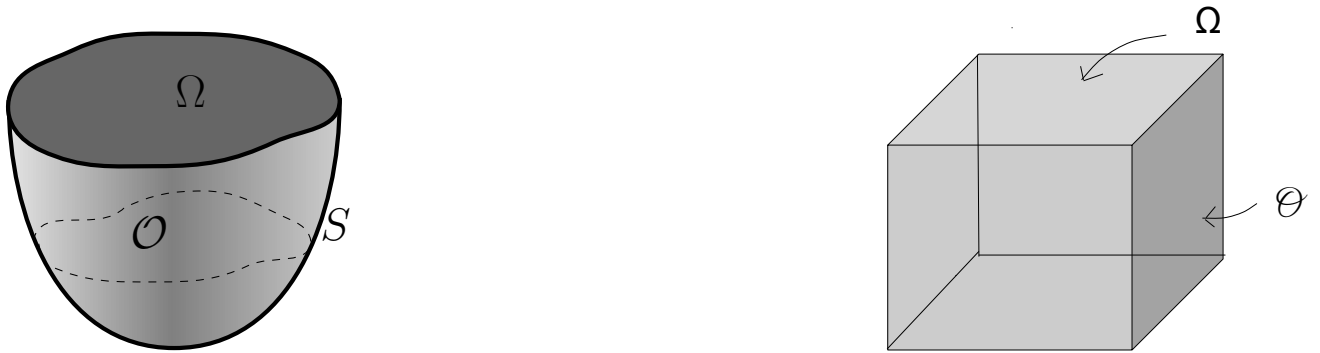
\begin{figure}[H]
\begin{subfigure}[h]{0.6\linewidth}
\begin{tikzpicture}[scale=1.1]
\draw[left color=black!10,right color=black!20,middle
color=black!50, ultra thick] (-2,0,0) to [out=0, in=180] (2,0,0)
to [out=270, in = 0] (0,-3,0) to [out=180, in =270] (-2,0,0);

\draw [fill=black!60, ultra thick] (-2,0,0) to [out=80,
in=205](-1.214,.607,0) to [out=25, in=180](0.2,.8,0) to [out=0,
in=155] (1.614,.507,0) to [out=335, in=100](2,0,0) to [out=270,
in=25] (1.214,-.507,0) to [out=205, in=0](-0.2,-.8,0) [out=180,
in=335] to (-1.614,-.607,0) to [out=155, in=260] (-2,0,0);

\draw [dashed, thin] (-1.7,-1.7,0) to [out=80, in=225](-.6,-1.3,0)
to [out=25, in=180](0.35,-1.1,0) to [out=0, in=155] (1.3,-1.4,0)
to [out=335, in=100](1.65,-1.7,0) to [out=270, in=25] (0.9,-2.0,0)
to [out=205, in=0](-0.2,-2.2,0) [out=180, in=335] to (-1.514,-2.0)
to [out=155, in=290] (-1.65,-1.7,0);

\node at (0.2,0.1,0) {{\LARGE$\Omega$}};

\node at (1.95,-1.5,0) {{\LARGE $S$}};

\node at (-0.3,-1.6,0) {{\LARGE $\mathcal{O}$}};
\end{tikzpicture}

\end{subfigure}
\hfill
\begin{subfigure}[h]{0.5\linewidth}
\includegraphics[width=\linewidth]{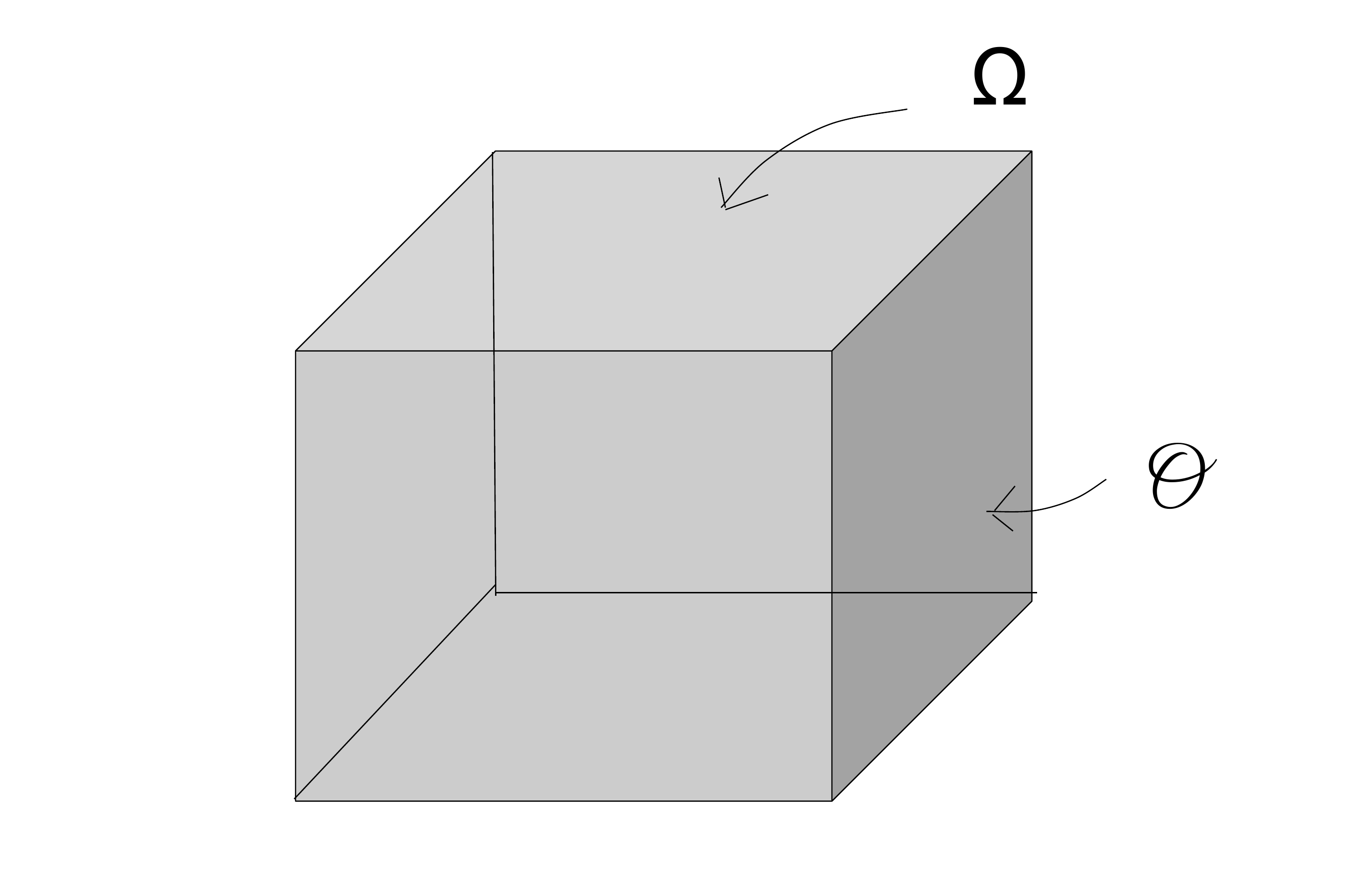}
\end{subfigure}%
\caption{Polyhedral Flow-Structure Geometries }
\end{figure}

The coupled PDE system which we will consider is the result of a
linearization which is undertaken in \cite{Chu2013-comp} and \cite{agw}:
Within the three-dimensional geometry $\mathcal{O}$, the compressible
Navier-Stokes equations are present, assuming the flow which they describe
to be barotropic. This system is linearized with respect to some reference
rest state of the form $\left\{ p_{\ast },\mathbf{U},\varrho _{\ast
}\right\} $: the pressure and density components ${p_{\ast },\varrho _{\ast }%
}$ are scalars, and the arbitrary ambient field $\mathbf{U}:\mathcal{O}%
\rightarrow \mathbb{R}^{3}$ 
\begin{equation*}
\mathbf{U}%
(x_{1},x_{2},x_{3})=[U_{1}(x_{1},x_{2},x_{3}),U_{2}(x_{1},x_{2},x_{3}),U_{3}(x_{1},x_{2},x_{3})]
\end{equation*}%
is given.

In \cite{Chu2013-comp} and \cite{agw}, we already see that non-critical lower order terms are
deleted, and the aforesaid pressure and density reference constants are set
equal to unity. Thus, we are presented with the following system of
equations, in solution variables $\mathbf{u}(x_{1},x_{2},x_{3},t)$ (flow
velocity), $p(x_{1},x_{2},x_{3},t)$ (pressure), $w(x_{1},x_{2},t)$ (elastic
plate displacement) and $w_{t}(x_{1},x_{2},t)$ (elastic plate velocity): 
\begin{align}
& \left\{ 
\begin{array}{l}
p_{t}+\mathbf{U}\cdot \nabla p+div~\mathbf{u}=0~\text{ in }~\mathcal{O}%
\times (0,\infty ) \\ 
\mathbf{u}_{t}+\mathbf{U}\cdot \nabla \mathbf{u}-div~\sigma (\mathbf{u}%
)+\eta \mathbf{u}+\nabla p=0~\text{ in }~\mathcal{O}\times (0,\infty ) \\ 
(\sigma (\mathbf{u})\mathbf{n}-p\mathbf{n})\cdot \boldsymbol{\tau }=0~\text{
on }~\partial \mathcal{O}\times (0,\infty ) \\ 
\mathbf{u}\cdot \mathbf{n}=0~\text{ on }~S\times (0,\infty ) \\ 
\mathbf{u}\cdot \mathbf{n}=w_{t}+\kappa \mathbf{U}\cdot \nabla w\text{ \ \
on }~\Omega \times (0,\infty )\text{ (where }\kappa =0\text{ or }\kappa =1%
\text{)}%
\end{array}%
\right.  \label{1} \\
& \left\{ 
\begin{array}{l}
w_{tt}+\Delta ^{2}w+\left[ 2\nu \partial _{x_{3}}(\mathbf{u})_{3}+\lambda 
\text{div}(\mathbf{u})-p\right] _{\Omega }=0~\text{ on }~\Omega \times
(0,\infty ) \\ 
w=\frac{\partial w}{\partial \nu }=0~\text{ on }~\partial \Omega \times
(0,\infty )%
\end{array}%
\right.  \label{2} \\
& 
\begin{array}{c}
\left[ p(0),\mathbf{u}(0),w(0),w_{t}(0)\right] =\left[ p_{0},\mathbf{u}%
_{0},w_{0},w_{1}\right]%
\end{array}
\label{3}
\end{align}%
(where above, $\nabla w=\left[ w_{x_{1}},w_{x_{2}},0\right] $; see \cite%
{material}). We note that the flow linearization is taken with respect to
some generally inhomogeneous compressible Navier-Stokes system; thus, $%
\mathbf{U}$ need not be divergence free generally; there are also initially
forcing terms in the pressure and flow equations (and energy level terms)
which we have neglected, since they do not effect the current analysis. This
flow-structure system is a generalization of that considered by the late
Igor Chueshov in \cite{Chu2013-comp} with therein fixed vector field $%
\mathbf{U}=\mathbf{0}$. In contrast, the PDE system (\ref{1})-(\ref{3})
depends upon a generally \emph{non-zero}, fixed, ambient vector field $%
\mathbf{U}$ about which the linearization takes place. The quantity $\eta >0$
represents a drag force of the domain on the viscous flow. In addition, the
quantity $\mathbf{\tau }$ in (\ref{1}) is in the space $TH^{1/2}(\partial 
\mathcal{O)}$ of tangential vector fields of Sobolev index 1/2; that is,%
\begin{equation}
\mathbf{\tau }\in TH^{1/2}(\partial \mathcal{O)=}\{\mathbf{v}\in \mathbf{H}^{%
\frac{1}{2}}(\partial \mathcal{O})~:~\left. \mathbf{v}\right\vert _{\partial 
\mathcal{O}}\cdot \mathbf{n}=0~\text{ on }~\partial \mathcal{O}\}.
\end{equation}
(See e.g., p.846 of \cite{buffa2}.) In addition, we take ambient field $%
\mathbf{U}$ to be in the space%
\begin{equation}
\mathbf{V}_{0}=\{\mathbf{v}\in \mathbf{H}^{1}(\mathcal{O})~:~\left. \mathbf{v%
}\right\vert _{\partial \mathcal{O}}\cdot \mathbf{n}=0~\text{ on }~\partial 
\mathcal{O}\}  \label{V_0}
\end{equation}%
(This vanishing of the boundary for ambient fields is a standard assumption
in compressible flow literature; see \cite{dV},\cite{valli},\cite{decay},%
\cite{spectral}.)

Moreover, the \textit{stress and strain tensors} in the flow PDE component
of (\ref{1})-(\ref{3}) are defined respectively as 
\begin{equation*}
\sigma (\mathbf{\mu })=2\nu \epsilon (\mathbf{\mu })+\lambda \lbrack
I_{3}\cdot \epsilon (\mathbf{\mu })]I_{3};\text{ \ }\epsilon _{ij}(\mathbf{%
\mu })=\dfrac{1}{2}\left( \frac{\partial \mathbf{\mu }_{j}}{\partial x_{i}}+%
\frac{\partial \mathbf{\mu }_{i}}{\partial x_{j}}\right) \text{, \ }1\leq
i,j\leq 3,
\end{equation*}%
where \textit{Lam\'{e} Coefficients }$\lambda \geq 0$ and $\nu >0$. The
associated finite energy space will be 
\begin{equation}
\mathcal{H}\equiv L^{2}(\mathcal{O})\times \mathbf{L}^{2}(\mathcal{O})\times
H_{0}^{2}(\Omega )\times L^{2}(\Omega )
\end{equation}%
which is a Hilbert space, topologized by the following standard inner
product: 
\begin{equation}
(\mathbf{y}_{1},\mathbf{y}_{2})_{\mathcal{H}}=(p_{1},p_{2})_{L^{2}(\mathcal{O%
})}+(\mathbf{u}_{1},\mathbf{u}_{2})_{\mathbf{L}^{2}(\mathcal{O})}+(\Delta
w_{1},\Delta w_{2})_{L^{2}(\Omega )}+(v_{1},v_{2})_{L^{2}(\Omega )}
\end{equation}%
for any $\mathbf{y}_{i}=(p_{i},\mathbf{u}_{i},w_{i},v_{i})\in \mathcal{H}%
,~i=1,2.$

\medskip

\begin{remark}
As we noted in \cite{agw}, the flow PDE boundary conditions which are in (%
\ref{1}) are the so-called \emph{impermeability}-slip conditions \cite%
{bolotin,chorin-marsden}: namely, no flow passes through the boundary -- in
particular, the normal component of the flow field $\mathbf{u}$ on the
active boundary portion $\Omega $ matches the plate velocity $w_{t}$ -- and
on $\partial \mathcal{O}$ there is no stress in the tangential direction $%
\tau $.
\end{remark}

\section{Functional Setting of the Problem}

\bigskip\ Throughout, for a given domain $D$, the norm of corresponding
space $L^{2}(D)$ will be denoted as $||\cdot ||_{D}$ (or simply $||\cdot ||$
when the context is clear). Inner products in $L^{2}(\mathcal{O})$ or $%
\mathbf{L}^{2}(\mathcal{O})$ will be denoted by $(\cdot ,\cdot )_{\mathcal{O}%
}$, whereas inner products $L^{2}(\partial \mathcal{O})$ will be written as $%
\langle \cdot ,\cdot \rangle _{\partial \mathcal{O}}$. We will also denote
pertinent duality pairings as $\left\langle \cdot ,\cdot \right\rangle
_{X\times X^{\prime }}$ for a given Hilbert space $X$. The space $H^{s}(D)$
will denote the Sobolev space of order $s$, defined on a domain $D$; $%
H_{0}^{s}(D)$ will denote the closure of $C_{0}^{\infty }(D)$ in the $%
H^{s}(D)$-norm $\Vert \cdot \Vert _{H^{s}(D)}$. We make use of the standard
notation for the boundary trace of functions defined on $\mathcal{O}$, which
are sufficently smooth: i.e., for a scalar function $\phi \in H^{s}(\mathcal{%
O})$, $\frac{1}{2}<s<\frac{3}{2}$, $\gamma (\phi )=\phi \big|_{\partial 
\mathcal{O}},$ which is a well-defined and surjective mapping on this range
of $s$, owing to the Sobolev Trace Theorem on Lipschitz domains (see e.g., 
\cite{necas}, or Theorem 3.38 of \cite{Mc}).

With respect to the above setting, the PDE system given in (\ref{1})-(\ref{3}%
) can be written as an ODE in Hilbert space $\mathcal{H}.$ That is, if $\Phi
(t)=\left[ p,\mathbf{u},w,w_{t}\right] \in C([0,T];\mathcal{H})$ solves the
problem (\ref{1})-(\ref{3}), then, for the respective cases $\kappa =0$ or $%
\kappa =1$ there is a modeling operator $\mathcal{A_{\kappa}}:D(\mathcal{A_{\kappa}})\subset 
\mathcal{H}\rightarrow \mathcal{H}$ such that $\Phi (\cdot )$ satisfies 
\begin{eqnarray}
\dfrac{d}{dt}\Phi (t) &=&\mathcal{A}_{\kappa }\Phi (t);  \notag \\
\Phi (0) &=&\Phi _{0}  \label{ODE}
\end{eqnarray}%
Here $\mathcal{A}_{\kappa }:D(\mathcal{A}_{\kappa })\subset \mathcal{H}%
\rightarrow \mathcal{H}$ is defined as follows:

\begin{equation}
\mathcal{A}_{\kappa }=\left[ 
\begin{array}{cccc}
-\mathbf{U}\mathbb{\cdot }\nabla (\cdot ) & -\text{div}(\cdot ) & 0 & 0 \\ 
-\mathbb{\nabla (\cdot )} & \text{div}\sigma (\cdot )-\eta I-\mathbf{U}%
\mathbb{\cdot \nabla (\cdot )} & 0 & 0 \\ 
0 & 0 & 0 & I \\ 
\left. \left[ \cdot \right] \right\vert _{\Omega } & -\left[ 2\nu \partial
_{x_{3}}(\cdot )_{3}+\lambda \text{div}(\cdot )\right] _{\Omega } & -\Delta
^{2} & 0%
\end{array}%
\right] ;  \label{AAA}
\end{equation}

\begin{equation*}
D(\mathcal{A}_{\kappa })=\{(p_{0},\mathbf{u}_{0},w_{1},w_{2})\in L^{2}(%
\mathcal{O})\times \mathbf{H}^{1}(\mathcal{O})\times H_{0}^{2}(\Omega
)\times H_{0}^{2}(\Omega )~:~\text{properties }(A.i)\text{--}(A.v)~~\text{%
hold}\},
\end{equation*}%
where

\begin{enumerate}
\item[(A.i)] $\mathbf{U}\cdot \nabla p_{0}\in L^{2}(\mathcal{O}).$

\item[(A.ii)] $\text{div}~\sigma (\mathbf{u}_{0})-\nabla p_{0}\in \mathbf{L}%
^{2}(\mathcal{O})$. (Consequently, we infer the boundary trace regularity 
\newline
$\left[ \sigma (\mathbf{u}_{0})\mathbf{n}-p_{0}\mathbf{n}\right] _{\partial 
\mathcal{O}}\in \mathbf{H}^{-\frac{1}{2}}(\partial \mathcal{O})$.)

\item[(A.iii)] $-\Delta ^{2}w_{1}-\left[ 2\nu \partial _{x_{3}}(\mathbf{u}%
_{0})_{3}+\lambda \text{div}(\mathbf{u}_{0})\right] _{\Omega }+\left.
p_{0}\right\vert _{\Omega }\in L^{2}(\Omega ).$

\item[(A.iv)] $\left( \sigma (\mathbf{u}_{0})\mathbf{n}-p_{0}\mathbf{n}%
\right) \bot ~TH^{1/2}(\partial \mathcal{O})$. That is, 
\begin{equation*}
\left\langle \sigma (\mathbf{u}_{0})\mathbf{n}-p_{0}\mathbf{n},\mathbf{\tau }%
\right\rangle _{\mathbf{H}^{-\frac{1}{2}}(\partial \mathcal{O})\times 
\mathbf{H}^{\frac{1}{2}}(\partial \mathcal{O})}=0\text{ \ for every }\mathbf{%
\tau }\in TH^{1/2}(\partial \mathcal{O})
\end{equation*}%
(and so $\left( \sigma (\mathbf{u}_{0})\mathbf{n}-p_{0}\mathbf{n}\right)
\cdot \mathbf{\tau }=0$ in the sense of distributions; see Remark 3.1 of 
\cite{agw}).

\item[(A.v)] The flow velocity component $\mathbf{u}_{0}=\mathbf{f}_{0}+%
\widetilde{\mathbf{f}}_{0}$, where $\mathbf{f}_{0}\in \mathbf{V}_{0}$ and $%
\widetilde{\mathbf{f}}_{0}\in \mathbf{H}^{1}(\mathcal{O})$ satisfies%
\footnote{%
The existence of an $\mathbf{H}^{1}(\mathcal{O})$-function $\widetilde{%
\mathbf{f}}_{0}$ with such a boundary trace on Lipschitz domain $\mathcal{O}$
is assured; see e.g., Theorem 3.33 of \cite{Mc}.}%
\begin{equation*}
\left. \widetilde{\mathbf{f}}_{0}\right\vert _{\partial \mathcal{O}}=%
\begin{cases}
0 & ~\text{ on }~S \\ 
\left( w_{2}+\kappa \mathbf{U}\cdot \nabla w_{1}\right) \mathbf{n} & ~\text{
on}~\Omega%
\end{cases}%
\end{equation*}%
\noindent (and so $\left. \mathbf{f}_{0}\right\vert _{\partial \mathcal{O}%
}\in TH^{1/2}(\partial \mathcal{O})$).
\end{enumerate}

\bigskip

In \cite{agw} and \cite{material}, it was shown that solutions to the
compressible flow-structure PDE system (\ref{1})-(\ref{3}), again for $%
\kappa =0$ or $\kappa =1$, with initial data in said finite energy space $%
\mathcal{H}$, can be associated with a strongly continuous semigroup $%
\left\{ e^{\mathcal{A}_{\kappa }t}\right\} \,_{t\geq 0}\subset \mathcal{L}(%
\mathcal{H)}$ which yields the following wellposedness result:

\medskip

\begin{theorem}
\label{well}  (See Theorem 3.1 of \cite{agw} and Theorem 5.1 of \cite{material}.) Assume that ambient vector field $\mathbf{U}%
\in \mathbf{V}_{0}\cap \mathbf{H}^{3}(\mathcal{O})$ (when $k=0$ and $k=1.$) Additionally, in the case $k=1,$ let $\mathbf{U}|_{\Omega}\in C^2(\bar{\Omega}).$ 

\textbf{(i)} (Wellposedness) The flow-structure
operator $\mathcal{A}_{\kappa }:D(\mathcal{A}_{\kappa })\subset \mathcal{H}%
\rightarrow \mathcal{H}$, generates a $C_{0}$-semigroup on $\mathcal{H}$.
Accordingly, the solution of (\ref{1})-(\ref{3}) -- with initial data $%
[p_{0},\mathbf{u}_{0},w_{0},w_{1}]\in $ $\mathcal{H}$ -- may be given by 
\begin{equation}
\begin{bmatrix}
p(t) \\ 
\mathbf{u}(t) \\ 
w(t) \\ 
w_{t}(t)%
\end{bmatrix}%
=e^{\mathcal{A}_{\kappa }t}%
\begin{bmatrix}
p_{0} \\ 
\mathbf{u}_{0} \\ 
w_{0} \\ 
w_{1}%
\end{bmatrix}%
\in C[(0,\infty );\mathcal{H}].  \label{w_1}
\end{equation}

\textbf{(ii)} (Energy relation) We define 
\begin{equation}
E(t)=\frac{1}{2}\left[ \left\Vert p(t)\right\Vert _{\mathcal{O}%
}^{2}+\left\Vert \mathbf{u}(t)\right\Vert _{\mathcal{O}}^{2}+\left\Vert
\Delta w(t)\right\Vert _{{\Omega }}^{2}+\left\Vert w_{t}(t)\right\Vert _{{%
\Omega }}^{2}\right] ,  \label{energy_0}
\end{equation}
where $\left[ p(t),\mathbf{u}(t),\Delta w(t),w_{t}(t)\right] $ is the
solution of compressible flow-structure system (\ref{1})-(\ref{3}), given
explicitly by (\ref{w_1}). Then for all $0\leq s<t$, we have the relation 
\begin{equation*}
E(t)+\int\limits_{s}^{t}\left[ (\sigma (\mathbf{u(}\tau \mathbf{)}),\epsilon
(\mathbf{u(}\tau \mathbf{))}_{\mathcal{O}}+\eta \left\Vert \mathbf{u(}\tau 
\mathbf{)}\right\Vert _{\mathcal{O}}^{2}\right] d\tau =E(s)+\frac{1}{2}%
\int\limits_{s}^{t}\int\limits_{\mathcal{O}}{div}(\mathbf{U})\left[
\left\vert p(\tau )\right\vert ^{2}+\left\vert \mathbf{u(}\tau \mathbf{)}%
\right\vert ^{2}\right] d\mathcal{O} d\tau  
\end{equation*}
\begin{equation}
+\kappa \int_{s}^{t}\left\langle \left[ 2\nu \partial _{x_{3}}(\mathbf{u}%
)_{3}+\lambda \text{div}(\mathbf{u})-p\right] _{\Omega },\mathbf{U}\cdot
\nabla w\right\rangle _{\Omega }d\tau. \label{energy-b}
\end{equation}
\end{theorem}

From the expression (\ref{energy-b}), it is seen that the generator $\mathcal{A%
}_{\kappa }:D(\mathcal{A}_{\kappa })\subset \mathcal{H}\rightarrow \mathcal{H%
}$ does not dissipate the energy of the system (\ref{1})-(\ref{3}). None the less, in \cite{pelin-george}, we establish that solutions of (\ref{1})-(\ref{3}) decay uniformly, with
respect to intial data which is $\mathcal{H}$-orthogonal to the one
dimensional null space of $\mathcal{A}_{0}$ (see Theorem 2 (ii) therein). 

The main intent of the present work is to \textbf{(i)} discern the null space for $%
\mathcal{A}_{1}:D(\mathcal{A}_{1})\subset \mathcal{H}\rightarrow \mathcal{H}$, under appropriate assumptions, and \textbf{(ii)} give an alternative proof of uniform
stability for the (material derivative-free) semigroup $\left\{ e^{\mathcal{A%
}_{0}t}\right\} _{t\geq 0}\subset \mathcal{L}(\mathcal{H})$, which in
contrast to the frequency domain approach of \cite{pelin-george}, is based
on the invocation of appropriate multipliers in the \emph{time domain}.

\bigskip

\section{The Null Space of $\mathcal{A}_{\protect\kappa }:D(\mathcal{A}_{%
\protect\kappa })\subset \mathcal{H}\rightarrow \mathcal{H}$}

In the process of analyzing the long term dynamics of the system (\ref{1})-(\ref{3}),  we need to avoid steady states so as to reasonably consider the possibility of finite energy solutions tending to the zero state at infinity. In fact, in the case $k=0$, the fact that zero is an eigenvalue for the generator $\mathcal{A}_{0}$ was proved in \cite{pelin-george}. Also an explicit characterization for the corresponding zero eigenspace was given there, again in the case $k=0.$ However, in the presence of the unbounded --material derivative-- term $\mathbf{U} \cdot \nabla w$ it is not at all clear a priori that the case $\kappa =1$ should give rise to the same zero eigenspace. In this section, we give a positive answer to this question: indeed, zero is also an eigenvalue for the generator $\mathcal{A}_{1}$ whose null space is identical with that for $\mathcal{A}_{0}.$ (As we said, this spectral information will be needed in our future work on longtime behavior properties of nonlinear compressible flow-structure PDE dynamics, with material derivative term in place.)  Before giving this result let us recall the following lemma (see Lemma 10 of \cite{pelin-george}):
\begin{lemma}
\label{SS} Let $\mathbf{U}\in \mathbf{V}_{0}\cap \mathbf{H}^{3}(\mathcal{O})$
and $\left\Vert {div}(\mathbf{U})\right\Vert _{\infty }\leq \mathfrak{C}%
_{0}^{\ast }(\Psi )$ (sufficiently small), where $\mathfrak{C}_{0}^{\ast
}(\Psi )$ is a positive constant with $\Psi =\Psi \left( \mathbf{U}\right)
\equiv \left[ \left\Vert \mathbf{U}\right\Vert _{\mathbf{H}^{3}(\mathcal{O}%
)}+1\right] $. Then one has the following:

The subspace $Null(\mathcal{A}_{0}\mathcal{)}\subset \mathcal{H}$ of the
flow-structure generator $\mathcal{A}_{0}:D(\mathcal{A}_{0})\subset \mathcal{%
H}\rightarrow \mathcal{H}$ is one dimensional. In particular, $Null(\mathcal{%
A}_{0}\mathcal{)}$ is given explicity as%
\begin{equation}
Null(\mathcal{A}_{0}\mathcal{)}=Span\left\{ \left[ 
\begin{array}{c}
1 \\ 
0 \\ 
{{\mathring{A}}^{-1}(1)} \\ 
0%
\end{array}%
\right] \right\} ,  \label{N1}
\end{equation}%
where $\mathring{A}:L^{2}(\Omega )\rightarrow L^{2}(\Omega )$ is the
elliptic operator%
\begin{equation}
\mathring{A}\varpi =\Delta ^{2}\varpi \text{, with }D(\mathring{A})=\{w\in
H_{0}^{2}(\Omega ):\Delta ^{2}w\in L^{2}(\Omega )\}.\text{\ \ }
\label{angst}
\end{equation}
\end{lemma}
\smallskip
Now, we give the same result for the generator $\mathcal{A}_{1}$ but we note that if we wish for the generator $\mathcal{A}_{1}:D(\mathcal{A}%
_{1})\subset \mathcal{H}\rightarrow \mathcal{H}$ of the flow-structure PDE
system to have the same one dimensional null space in (\ref{N1}), as for $%
\mathcal{A}_{0}$, despite the additional (unbounded) material derivative
term, then: from the domain criterion (A.v) and (\ref{N1}), one must
necessarily have, for $\kappa =1$, the relation%
\begin{equation}
\mathbf{U}\cdot \nabla {{\mathring{A}}^{-1}(1)=0}\text{  on  }\Omega
\label{cc}
\end{equation}%
where the biharmonic operator $\mathring{A}$ is as in (\ref{angst}). That
is, $$\left[ U_{1}(x_{1},x_{2},0),U_{2}(x_{1},x_{2},0\right] \cdot \left[ 
\frac{\partial }{\partial x_{1}}{{\mathring{A}}^{-1}(1),}\frac{\partial }{%
\partial x_{2}}{{\mathring{A}}^{-1}(1)}\right] =0.$$
\smallskip

\begin{lemma}
\label{SS_2}(a) Suppose that the ambient field $\mathbf{U}$ satisfies (\ref%
{cc}), and moreover $\left\{ \left\Vert {div}(\mathbf{U})\right\Vert
_{\infty },\left\Vert \left. \mathbf{U}\right\vert _{\Omega }\right\Vert _{%
\boldsymbol{H}^{\frac{1}{2}}(\Omega )\cap \mathbf{L}^{\infty }(\Omega
)}\right\} \leq \mathfrak{C}_{0}^{\ast }(\Psi )$ (sufficiently small), where 
$\mathfrak{C}_{0}^{\ast }(\Psi )$ is a positive constant with $\Psi =\Psi
\left( \mathbf{U}\right) \equiv \left[ \left\Vert \mathbf{U}\right\Vert _{%
\mathbf{H}^{3}(\mathcal{O})}+1\right] $. Then, as in Lemma \ref{SS}, the
subspace $Null(\mathcal{A}_{1}\mathcal{)}\subset \mathcal{H}$ of the
(when the material derivative is taken into account) flow-structure generator $\mathcal{A}_{1}$ is given by 
\begin{equation}
Null(\mathcal{A}_{1})\mathcal{=}Null(\mathcal{A}_{0}\mathcal{)}=Span\left\{ %
\left[ 
\begin{array}{c}
1 \\ 
0 \\ 
{{\mathring{A}}^{-1}(1)} \\ 
0%
\end{array}%
\right] \right\}   \label{N1.5}
\end{equation}

(b) The orthogonal complement of $Null(\mathcal{A}_{\kappa }\mathcal{)}$ ($%
\kappa =0$ or $\kappa =1$) admits of the characterization 
\begin{equation}
\lbrack Null(\mathcal{A}_{\kappa }\mathcal{)}]^{\bot }=\mathcal{\{}[p_{0},%
\mathbf{u}_{0},w_{1},w_{2}]\in \mathcal{H}:\int\limits_{\mathcal{O}}p_{0}d%
\mathcal{O}+\int\limits_{\Omega }w_{1}d\Omega =0\mathcal{\}}.  \label{N2}
\end{equation}
\end{lemma}

\noindent \textbf{Proof of Lemmas \ref{SS} and \ref{SS_2}} \newline\newline
\ We focus here on the touchier case $\kappa =1$, inasmuch as the details
for $\kappa =0$ are subsumed in the material derivative case. (See also \cite%
{pelin-george}.)

\medskip

\noindent Suppose $\Phi =[p_{0},\mathbf{u}_{0},w_{1},w_{2}]\in D(\mathcal{A}_{1})$ is
a solution of 
\begin{equation}
\mathcal{A}_{1}\Phi =\mathbf{0}\text{,}  \label{zero}
\end{equation}%
where $\mathcal{A}_{1}:D(\mathcal{A}_{1})\subset \mathcal{H}\rightarrow 
\mathcal{H}$ is as given in (\ref{AAA}). (Without loss of generality, we
take $\Phi $ to be real-valued.) With respect to the pressure component, we
invoke the $L^{2}$-decomposition 
\begin{equation}
p_{0}=q_{0}+c_{0}\text{,}  \label{nu_1.5}
\end{equation}
where 
\begin{equation}
q_{0}\text{ satisfies }\int\limits_{\mathcal{O}}q_{0}d\mathcal{O}=0\text{, \
and \ }c_{0}=\text{constant.}  \label{nu_2}
\end{equation}

\medskip

\noindent Therewith, in PDE terms, the abstract relation then becomes

\begin{eqnarray}
&&\left\{ 
\begin{array}{l}
-\mathbf{U}\cdot \nabla q_{0}-div~\mathbf{u}_{0}=0\text{ \ in }~\mathcal{O}
\\ 
-\nabla q_{0}+div~\sigma (\mathbf{u}_{0})-\eta \mathbf{u}_{0}-\mathbf{U}%
\cdot \nabla \mathbf{u}_{0}=0~\text{ \ in }~\mathcal{O} \\ 
w_{2}=0\text{ \ on }~\Omega \\ 
-\Delta ^{2}w_{1}-\left[ 2\nu \partial _{x_{3}}(\mathbf{u}_{0})_{3}+\lambda 
\text{div}(\mathbf{u}_{0})-p_{0}\right] _{\Omega }=0~\text{ on }~\Omega%
\end{array}%
\right.  \label{eq_1} \\
&&  \notag \\
&&\left\{ 
\begin{array}{l}
(\sigma (\mathbf{u}_{0})\mathbf{n}-p_{0}\mathbf{n})\cdot \boldsymbol{\tau }=0%
\text{ \ on }~\partial \mathcal{O} \\ 
\mathbf{u}_{0}\cdot \mathbf{n}=0~\text{ on }~S \\ 
\mathbf{u}_{0}\cdot \mathbf{n}=\mathbf{U}\cdot \nabla w_{1}~\text{ on }%
~\Omega \\ 
w_{1}=\frac{\partial w_{1}}{\partial \nu }=0~\text{ on }~\partial \Omega .%
\end{array}%
\right.  \label{eq_2}
\end{eqnarray}

\medskip \noindent We have immediately then,%
\begin{equation}
w_{2}=0\text{.}  \label{nu_0}
\end{equation}%
\noindent Secondly, we multiply the pressure PDE component in (\ref{eq_1})-(%
\ref{eq_2}) by $p_{0}$ and the fluid PDE component by $\mathbf{u}%
_{0}$. Subsequent integrations and integrations by parts, and a
consideration of domain criteria (A.iv), (A.v) (and (\ref{nu_0})), yield then%
\begin{equation}
\begin{array}{l}
\left( \sigma (\mathbf{u}_{0}),\epsilon (\mathbf{u}_{0})\right) _{\mathcal{O}%
}+\eta \left\Vert \mathbf{u}_{0}\right\Vert _{\mathcal{O}}^{2}=\frac{1}{2}%
\left( {div}(\mathbf{U})p_{0},p_{0}\right) _{\mathcal{O}}+\frac{1}{2}\left( {%
div}(\mathbf{U})\mathbf{u}_{0},\mathbf{u}_{0}\right) _{\mathcal{O}} \\ 
\text{ \ \ \ }\ -\left\langle c_{0},\mathbf{U}\cdot \nabla
w_{1}\right\rangle _{\Omega }+\left\langle \left[ 2\nu \partial _{x_{3}}(%
\mathbf{u}_{0})_{3}+\lambda \text{div}(\mathbf{u}_{0})-q_{0}\right] _{\Omega
},u_{3}\right\rangle _{\Omega }%
\end{array}
\label{c1}
\end{equation}%
(For the last term on RHS, we are also using the fact that $\mathbf{n}%
=[0,0,1]$ on $\Omega$.) Therewith, combining the decomposition (\ref{nu_1.5}) with
Green's formula and the fact that $\mathbf{U}\in \mathbf{V}_{0}$, we then
obtain%
\begin{eqnarray}
&&\left( \sigma (\mathbf{u}_{0}),\epsilon (\mathbf{u}_{0})\right) _{\mathcal{%
O}}+\eta \left\Vert \mathbf{u}_{0}\right\Vert _{\mathcal{O}}^{2}=\left( {div}%
(\mathbf{U})c_{0},q_{0}\right) _{\mathcal{O}}  \notag \\
&&\text{ \ }-\left\langle c_{0},\mathbf{U}\cdot \nabla w_{1}\right\rangle
_{\Omega }+\left\langle \left[ 2\nu \partial _{x_{3}}(\mathbf{u}%
_{0})_{3}+\lambda \text{div}(\mathbf{u}_{0})-q_{0}\right] _{\Omega
},u_{3}\right\rangle _{\Omega }  \notag \\
&&\text{ \ }+\frac{1}{2}\left( {div}(\mathbf{U})q_{0},q_{0}\right) _{%
\mathcal{O}}+\frac{1}{2}\left( {div}(\mathbf{U})\mathbf{u}_{0},\mathbf{u}%
_{0}\right) _{\mathcal{O}}.  \label{p1}
\end{eqnarray}

\medskip \noindent To estimate the first term on right hand side of (\ref{p1}%
), we multiply the pressure equation in (\ref{eq_1}) by constant component $%
c_{0}$ of (\ref{nu_1.5}) and integrate over $\mathcal{O}$. This gives, 
\begin{equation*}
-\left( \mathbf{U}\cdot \nabla q_{0},c_{0}\right) =\left( {div}\mathbf{u%
}_{0},c_{0}\right) .
\end{equation*}%
Subsequently we integrate by parts both sides of this relation, while
bearing in mind the domain criterion (A.v) and (\ref{nu_0}) (and the fact
that $\mathbf{U}\in \mathbf{V}_{0}$), so as to have 
\begin{equation*}
\left( {{div}}(\mathbf{U})q_{0},c_{0}\right) _{\mathcal{O}%
}=\left\langle c_{0},\mathbf{U}\cdot \nabla w_{1}\right\rangle _{\Omega }.
\end{equation*}%
Applying this relation to (\ref{p1}) gives now%
\begin{eqnarray*}
\left( \sigma (\mathbf{u}_{0}),\epsilon (\mathbf{u}_{0})\right) _{\mathcal{O}%
}+\eta \left\Vert \mathbf{u}_{0}\right\Vert _{\mathcal{O}}^{2}
&=&\left\langle \left[ 2\nu \partial _{x_{3}}(\mathbf{u}_{0})_{3}+\lambda 
\text{div}(\mathbf{u}_{0})-q_{0}\right] _{\Omega },u_{3}\right\rangle
_{\Omega } \\
&&+\frac{1}{2}\left( {div}(\mathbf{U})q_{0},q_{0}\right) _{\mathcal{O}}+%
\frac{1}{2}\left( {div}(\mathbf{U})\mathbf{u}_{0},\mathbf{u}_{0}\right) _{%
\mathcal{O}};
\end{eqnarray*}%
and subsequently invoking the domain criterion (A.v) (for $\kappa =1$) and
the mechanical equation on $\Omega $, we have now%
\begin{eqnarray*}
\left( \sigma (\mathbf{u}_{0}),\epsilon (\mathbf{u}_{0})\right) _{\mathcal{O}%
}+\eta \left\Vert \mathbf{u}_{0}\right\Vert _{\mathcal{O}}^{2}
&=&\left\langle \left[ 2\nu \partial _{x_{3}}(\mathbf{u}_{0})_{3}+\lambda 
\text{div}(\mathbf{u}_{0})-q_{0}\right] _{\Omega },\mathbf{U}\cdot \nabla
\left( {{\mathring{A}}^{-1}(c}_{0}-\left[ 2\nu \partial _{x_{3}}(\mathbf{u}%
_{0})_{3}+\lambda \text{div}(\mathbf{u}_{0})-q_{0}\right] _{\Omega }{)}%
\right) \right\rangle _{\Omega } \\
&&+\frac{1}{2}\left( {div}(\mathbf{U})q_{0},q_{0}\right) _{\mathcal{O}}+%
\frac{1}{2}\left( {div}(\mathbf{U})\mathbf{u}_{0},\mathbf{u}_{0}\right) _{%
\mathcal{O}};
\end{eqnarray*}%
and subsequently invoking (\ref{cc}) yields then%
\begin{eqnarray}
\left( \sigma (\mathbf{u}_{0}),\epsilon (\mathbf{u}_{0})\right) _{\mathcal{O}%
}+\eta \left\Vert \mathbf{u}_{0}\right\Vert _{\mathcal{O}}^{2}
&=&-\left\langle \left[ 2\nu \partial _{x_{3}}(\mathbf{u}_{0})_{3}+\lambda 
\text{div}(\mathbf{u}_{0})-q_{0}\right] _{\Omega },\mathbf{U}\cdot \nabla
\left( {{\mathring{A}}^{-1}}\left[ 2\nu \partial _{x_{3}}(\mathbf{u}%
_{0})_{3}+\lambda \text{div}(\mathbf{u}_{0})-q_{0}\right] _{\Omega }\right)
\right\rangle _{\Omega }  \notag \\
&&+\frac{1}{2}\left( {div}(\mathbf{U})q_{0},q_{0}\right) _{\mathcal{O}}+%
\frac{1}{2}\left( {div}(\mathbf{U})\mathbf{u}_{0},\mathbf{u}_{0}\right) _{%
\mathcal{O}}.  \label{p2}
\end{eqnarray}

\bigskip \noindent To make use of this relation, we appeal to the estimate
provided by Theorem 2.4 or Remark 2.5 of \cite{temam}, for solution pair $%
\left( \mathbf{u}_{0},q_{0}\right) $ of (Stokes system) (\ref{eq_1}).
Namely, we have 
\begin{eqnarray*}
&&\left\Vert q_{0}\right\Vert _{\mathcal{O}}^{2}+\left( \sigma (\mathbf{u}%
_{0}),\epsilon (\mathbf{u}_{0})\right) _{\mathcal{O}}+\eta \left\Vert 
\mathbf{u}_{0}\right\Vert _{\mathcal{O}}^{2} \\
&\leq &C\left( \left\Vert \mathbf{U}\cdot \nabla \mathbf{u}_{0}\right\Vert _{%
\mathcal{O}}^{2}+\left\Vert {div}(\mathbf{u}_{0})\right\Vert _{\mathcal{O}%
}^{2}+\left\Vert \mathbf{u}_{0}\right\Vert _{\mathbf{H}^{\frac{1}{2}%
}(\partial \mathcal{O})}^{2}\right) 
\end{eqnarray*}%
(We are also implicitly using Korn's Inequality.) Estimating right hand side
by Korn's Inequality and the Sobolev Trace Theorem gives now%
\begin{eqnarray}
&&\left\Vert q_{0}\right\Vert _{\mathcal{O}}^{2}+\left( \sigma (\mathbf{u}%
_{0}),\epsilon (\mathbf{u}_{0})\right) _{\mathcal{O}}+\eta \left\Vert 
\mathbf{u}_{0}\right\Vert _{\mathcal{O}}^{2}  \notag \\
&\leq &\mathfrak{C}_{1}\Psi (\mathbf{U})\left[ \left( \sigma (\mathbf{u}%
_{0}),\epsilon (\mathbf{u}_{0})\right) _{\mathcal{O}}+\eta \left\Vert 
\mathbf{u}_{0}\right\Vert _{\mathcal{O}}^{2}\right] ,  \label{p22}
\end{eqnarray}%
where the term $\Psi (\mathbf{U})$ is as given in the statement of Lemma \ref%
{SS_2}, and positive constant $\mathfrak{C}_{1}$ is independent of $\mathbf{U%
}.$ (We are implicitly using here the Sobolev Imbedding Theorem with respect
to field $\mathbf{U}\in \mathbf{H}^{3}(\mathcal{O})$). Applying 
(\ref{p22}) to the right hand side of (\ref{p2}) gives then, 
\begin{equation*}
\left( \sigma (\mathbf{u}_{0}),\epsilon (\mathbf{u}_{0})\right) _{\mathcal{%
O}}+\eta \left\Vert \mathbf{u}_{0}\right\Vert _{\mathcal{O}}^{2} \leq
\end{equation*}
\begin{equation*}
 \left\vert \left\langle \left[ 2\nu \partial _{x_{3}}(\mathbf{u}%
_{0})_{3}+\lambda \text{div}(\mathbf{u}_{0})-q_{0}\right] _{\Omega },\mathbf{%
U}\cdot \nabla \left( {{\mathring{A}}^{-1}}\left[ 2\nu \partial _{x_{3}}(%
\mathbf{u}_{0})_{3}+\lambda \text{div}(\mathbf{u}_{0})-q_{0}\right] _{\Omega
}\right) \right\rangle _{\Omega }\right\vert  
\end{equation*}
\begin{equation*}
+{\mathfrak{C}_{2}(\Psi (\mathbf{U}))}\left\Vert {div}(\mathbf{U}%
)\right\Vert _{\infty }\left[ \left( \sigma (\mathbf{u}_{0}),\epsilon (%
\mathbf{u}_{0})\right) _{\mathcal{O}}+\eta \left\Vert \mathbf{u}%
_{0}\right\Vert _{\mathcal{O}}^{2}\right],
\end{equation*}%
whence we obtain%
\begin{eqnarray}
&&\left( \sigma (\mathbf{u}_{0}),\epsilon (\mathbf{u}_{0})\right) _{\mathcal{%
O}}+\eta \left\Vert \mathbf{u}_{0}\right\Vert _{\mathcal{O}}^{2}  \notag \\
&\leq &C\left\Vert \left. \mathbf{U}\right\vert _{\Omega }\right\Vert _{%
\boldsymbol{H}^{\frac{1}{2}}(\Omega )\cap \mathbf{L}^{\infty }(\Omega
)}\left\Vert \left[ 2\nu \partial _{x_{3}}(\mathbf{u}_{0})_{3}+\lambda \text{%
div}(\mathbf{u}_{0})-q_{0}\right] _{\Omega }\right\Vert _{H^{-\frac{1}{2}%
}(\Omega )}^{2}  \notag \\
&&+{\mathfrak{C}_{2}(\Psi (\mathbf{U}))}\left\Vert {div}(\mathbf{U}%
)\right\Vert _{\infty }\left[ \left( \sigma (\mathbf{u}_{0}),\epsilon (%
\mathbf{u}_{0})\right) _{\mathcal{O}}+\eta \left\Vert \mathbf{u}%
_{0}\right\Vert _{\mathcal{O}}^{2}\right] .  \label{p3.5}
\end{eqnarray}

\medskip

\noindent In turn, an integration by parts, the estimates (\ref{p22}), and (\ref{p3.5})
(and a rescaling of parameter $\epsilon >0$) give%
\begin{equation}
\left\Vert \sigma (\mathbf{u}_{0})\cdot \mathbf{n}-q_{0}\mathbf{n}%
\right\Vert _{\mathbf{H}^{-\frac{1}{2}}(\partial \mathcal{O})}^{2}\leq
C\left( \left\Vert q_{0}\right\Vert _{\mathcal{O}}^{2}+\left( \sigma (%
\mathbf{u}_{0}),\epsilon (\mathbf{u}_{0})\right) _{\mathcal{O}}+\eta
\left\Vert \mathbf{u}_{0}\right\Vert _{\mathcal{O}}^{2}+\left\Vert \mathbf{U}%
\cdot \nabla \mathbf{u}_{0}\right\Vert _{\mathcal{O}}^{2}\right) .
\label{p3.6}
\end{equation}%
Combining (\ref{p3.5}) and (\ref{p3.6}) gives now%
\begin{eqnarray}
&&\left\Vert \sigma (\mathbf{u}_{0})\cdot \mathbf{n}-q_{0}\mathbf{n}%
\right\Vert _{\mathbf{H}^{-\frac{1}{2}}(\partial \mathcal{O})}^{2}+\left(
\sigma (\mathbf{u}_{0}),\epsilon (\mathbf{u}_{0})\right) _{\mathcal{O}}+\eta
\left\Vert \mathbf{u}_{0}\right\Vert _{\mathcal{O}}^{2}  \notag \\
&\leq &{\mathfrak{C}_{3}(\Psi (\mathbf{U}))}\left( \left\Vert \left. \mathbf{%
U}\right\vert _{\Omega }\right\Vert _{\boldsymbol{H}^{\frac{1}{2}}(\Omega
)\cap \mathbf{L}^{\infty }(\Omega )}\left\Vert \left[ 2\nu \partial _{x_{3}}(%
\mathbf{u}_{0})_{3}+\lambda \text{div}(\mathbf{u}_{0})-q_{0}\right] _{\Omega
}\right\Vert _{H^{-\frac{1}{2}}(\Omega )}^{2}\right.   \notag \\
&&\left. +\left\Vert {div}(\mathbf{U})\right\Vert _{\infty }\left[ \left(
\sigma (\mathbf{u}_{0}),\epsilon (\mathbf{u}_{0})\right) _{\mathcal{O}}+\eta
\left\Vert \mathbf{u}_{0}\right\Vert _{\mathcal{O}}^{2}\right] \right) .
\label{p3.7}
\end{eqnarray}

\medskip

\noindent If we now specify%
\begin{equation}
\max \left\{ \left\Vert {div}(\mathbf{U})\right\Vert _{\infty },\left\Vert
\left. \mathbf{U}\right\vert _{\Omega }\right\Vert _{\boldsymbol{H}^{\frac{1%
}{2}}(\Omega )\cap \mathbf{L}^{\infty }(\Omega )}\right\} <\frac{1}{{%
\mathfrak{C}_{3}(\Psi (\mathbf{U}))}}\equiv {\mathfrak{C}_{0}^{\ast }(\Psi (%
\mathbf{U}));}  \label{p3.8}
\end{equation}%
where again $\Psi \left( \mathbf{U}\right) \equiv \left[ \left\Vert \mathbf{U%
}\right\Vert _{\mathbf{H}^{3}(\mathcal{O})}+1\right] $, we then deduce that

\begin{equation}
\mathbf{u}_{0}=\mathbf{0}\text{  and  }q_{0}=0\text{.}  \label{p4}
\end{equation}%
(We are also implicitly using Korn's Inequality). Consequently, via (\ref%
{nu_1.5}) we obtain 
\begin{equation}
p_{0}=c_{0}\text{, where }c_{0}\text{ is constant.}  \label{p5}
\end{equation}%
Returning to the mechanical equation on $\Omega $ in (\ref{eq_1}), we then
have 
\begin{equation}
w_{1}={{\mathring{A}}^{-1}(c_{0}),}  \label{p6}
\end{equation}%
where again $\mathring{A}^{-1}:L^{2}(\Omega )\rightarrow L^{2}(\Omega )$ is
the operator defined in (\ref{angst}). The relations (\ref{nu_0}), (\ref{p4}%
), (\ref{p5}) and (\ref{p6}) give the conclusion of (a) of Lemma \ref{SS_2} for $\kappa =1$. Given the definition of the $\mathcal{H}$-inner product as well as the
definition in (\ref{angst}) of $\mathring{A}:L^{2}(\Omega )\rightarrow
L^{2}(\Omega )$, the relation in (b) of Lemma \ref{SS_2} is immediate. \ \ \ $\square $

\medskip

\begin{remark}
Suppose the curvilinear polyhedron $\mathcal{O}$ is strictly an
\textquotedblleft edge domain\textquotedblright\ -- i.e., the geometry has
no corners -- such that \text{each point on an edge of }$\partial \mathcal{O}
$ \text{ is diffeomorphic to a wedge with opening} $<\pi /2$. Then an
example of an ambient field $\mathbf{U\in V}_{0}\cap \mathbf{H}^{3}(\mathcal{%
O})$ which meets the compatibility condition (\ref{cc})
can be constructed  as follows: Let boundary data $\mathbf{g}$ satisfy%
\begin{equation*}
\mathbf{g}=\left\{ 
\begin{array}{l}
\left[ 
\begin{array}{c}
-\frac{\partial }{\partial x_{2}}{{\mathring{A}}^{-1}(1)} \\ 
\frac{\partial }{\partial x_{1}}{{\mathring{A}}^{-1}(1)} \\ 
0%
\end{array}%
\right] \in \mathbf{H}^{s}(\Omega )\text{ \ on }\Omega  \\ 
\\ 
\mathbf{0}\text{ \ on }\partial \mathcal{O}.%
\end{array}%
\right. 
\end{equation*}%
Therewith, $\mathbf{g}\in \mathbf{H}^{3}(\partial \mathcal{O})$; see e.g.,
Theorem 3.33, p. 95, of \cite{Mc}. Here, we are implicitly using the fact that 
$\left. \frac{\partial }{\partial x_{i}}{{\mathring{A}}^{-1}(1)}\right\vert
_{\partial \Omega }=0$. Subsequently, we can invoke the well-established
elliptic regularity results on edge domains -- see; e.g., \cite{Dauge_3} --
for Neumann problem, with boundary data $\mathbf{g}$, so as to secure a $%
\mathbf{U}$ which meets the requirements of Lemma \ref{SS_2}, after an
appropriate scaling.
\end{remark}

\section{Exponential Stability -- A Time Domain Approach}

This section is devoted to giving an alternative proof for the exponential stability of the solutions to material derivative-free system (\ref{1})-(\ref{3}) (i.e., the case $k=0$.) Unlike the frequency domain approach followed in \cite{pelin-george} to obtain this stability result, we give the proof of Theorem (\ref%
{uniform}) in time domain, an approach which involves a gradient type multiplier by way of obtaining necessary energy estimates.
As we mentioned before, the alternative proof which we provide here can in principle be used in any long term analysis for corresponding nonlinear systems. In fact, the time domain approach which we develop here will be our key departure point in the forthcoming paper related to compact global attractors of the system (\ref{1})-(\ref{3}), in the presence of the material derivative, as well as given structural (von Karman) nonlinearity.  

With the notation adopted above, we
give the following exponential decay result for all initial data taken
from $[Null(\mathcal{A)}]^{\bot }$, which is a one dimensional subspace of $%
\mathcal{H}$ (see Lemma 6(b)):

\begin{theorem}
\label{uniform} Let the ambient vector field $\mathbf{U}\in \mathbf{V}%
_{0}\cap \mathbf{H}^{3}(\mathcal{O})$ and the geometrical assumptions in Condition \ref{cond} hold. Additionally, assume that $\left\Vert {div}(\mathbf{U}%
)\right\Vert _{\infty }$ is sufficiently small. Then the (material
derivative -free) $C_{0}$-semigroup $\left\{ \left. e^{\mathcal{A}%
_{0}t}\right\vert _{[Null(\mathcal{A}_{0}\mathcal{)}]^{\bot }}\right\}
_{t\geq 0}\in \mathcal{L}\left( [Null(\mathcal{A}_{0}\mathcal{)}]^{\bot
}\right) $ decays exponentially. In particular, there exist constants $%
\mathfrak{M}>0$ and $\omega >0$ such that the solution to the flow-structure
PDE system (\ref{1})-(\ref{3}), with initial data $\Phi _{0}=\left[ p_{0},%
\mathbf{u}_{0},w_{0}w_{1}\right] \in \lbrack Null(\mathcal{A}_{0}\mathcal{)}%
]^{\bot }$, obeys the uniform decay rate 
\begin{equation}
\left\Vert e^{\mathcal{A}_{0}t}\Phi _{0}\right\Vert _{\mathcal{H}}\leq 
\mathfrak{M}e^{-\omega t}\left\Vert \Phi _{0}\right\Vert _{\mathcal{H}}\text{%
, \ for all }t>0\text{.}  \label{exp}
\end{equation}
\end{theorem}

\begin{remark}
In our time domain approach, the geometrical assumptions on the polyhedral flow domain $\mathcal{O}$ are necessary, since in the course of proof of Theorem 8, we appeal to higher regularity results for the Neumann problem on domains with corners (see \cite{JK} and \cite{Dauge_3} for further details.)  They are analogous to the geometrical assumptions made in \cite{pelin-george}, in which the frequency domain approach requires the invocation of higher regularity results for static Stokes flow on corner domains; see \cite{dauge}.
\end{remark}

\begin{remark}
\label{inv}We note that for $\kappa =0$, one has the invariance of $\mathcal{%
A}_{0}$ and its associated $C_{0}$-semigroup $\left\{ e^{\mathcal{A}%
_{0}t}\right\} _{t\geq 0}$ on $[Null(\mathcal{A}_{0}\mathcal{)}]^{\bot }$;
in particular, $\left\{ \left. e^{\mathcal{A}_{0}t}\right\vert _{[Null(%
\mathcal{A)}]^{\bot }}\right\} \subset \mathcal{L}\left( [Null(\mathcal{A}%
_{0}\mathcal{)}]^{\bot }\right) $, for $\kappa =0.$ (See Proposition 13 of  
\cite{pelin-george}.)
\end{remark}
\medskip 
\textbf{Proof of Theorem \ref{uniform}}. \\

The proof is based on the application of a multiplier which is  intrinsic to the compressible flow-structure PDE system under consideration. In order to deal with the lack of $H^1$-regularity of the pressure variable $p$ -- which is partly manifested in the unbounded term $\mathbf{U}\cdot \nabla p$ -- this special multiplier exploits the compatibility condition given in (\ref{N2}) for any data of $\lbrack Null(\mathcal{A}_{0}\mathcal{)}%
]^{\bot }$ so as to ultimately enable us to obtain  the necessary observability estimate for the energy of the system.

We will consider the following
system (in the case $\kappa =0$ in (\ref{1})-(\ref{3}) ) and initial data $\Phi (0)=\left[ p(0),\mathbf{u}%
(0),w(0),w_{t}(0)\right] =\left[ p_{0},\mathbf{u}_{0},w_{0},w_{1}\right]
=\Phi _{0}\in \lbrack Null(\mathcal{A}_{0}\mathcal{)}]^{\bot }:$%
\begin{align}
& \left\{ 
\begin{array}{l}
p_{t}+\mathbf{U}\cdot \nabla p+div~\mathbf{u}=0~\text{ in }~\mathcal{O}%
\times (0,\infty ) \\ 
\mathbf{u}_{t}+\mathbf{U}\cdot \nabla \mathbf{u}-div~\sigma (\mathbf{u}%
)+\eta \mathbf{u}+\nabla p=0~\text{ in }~\mathcal{O}\times (0,\infty ) \\ 
(\sigma (\mathbf{u})\mathbf{n}-p\mathbf{n})\cdot \boldsymbol{\tau }=0~\text{
on }~\partial \mathcal{O}\times (0,\infty )\text{ \ \ for all }\boldsymbol{%
\tau }\in TH^{1/2}(\partial \mathcal{O}) \\ 
\mathbf{u}\cdot \mathbf{n}=0~\text{ on }~S\times (0,\infty ) \\ 
\mathbf{u}\cdot \mathbf{n}=w_{t}\text{ \ \ on }%
~\Omega \times (0,\infty )%
\end{array}%
\right.   \label{k1} \\
& \left\{ 
\begin{array}{l}
w_{tt}+\Delta ^{2}w+\left[ 2\nu \partial _{x_{3}}(\mathbf{u})_{3}+\lambda 
\text{div}(\mathbf{u})-p\right] _{\Omega }=0~\text{ on }~\Omega \times
(0,\infty ) \\ 
w=\frac{\partial w}{\partial \nu }=0~\text{ on }~\partial \Omega \times
(0,\infty )%
\end{array}%
\right.   \label{k2} \\
& 
\begin{array}{c}
\left[ p(0),\mathbf{u}(0),w(0),w_{t}(0)\right] =\left[ p_{0},\mathbf{u}%
_{0},w_{0},w_{1}\right] 
\end{array}%
\in \lbrack Null(\mathcal{A}_{0}\mathcal{)}]^{\bot }.  \label{k3}
\end{align}%
From Theorem \ref{well}(ii), we have the following energy relation:%
\begin{equation}
E(t)+\int\limits_{s}^{t}\left[ (\sigma (\mathbf{u(}\tau \mathbf{)}),\epsilon
(\mathbf{u(}\tau \mathbf{))}_{\mathcal{O}}+\eta \left\Vert \mathbf{u(}\tau 
\mathbf{)}\right\Vert _{\mathcal{O}}^{2}\right] d\tau =E(s)+\frac{1}{2}%
\int\limits_{s}^{t}\int\limits_{\mathcal{O}}{div}(\mathbf{U})\left[
\left\vert p(\mathbf{\tau })\right\vert ^{2}+\left\vert \mathbf{u(\tau )}%
\right\vert ^{2}\right] d\mathcal{O} d\tau   \label{E2}
\end{equation}
where $E(t)$ is as given in (\ref{energy_0}).

\medskip

By the classic \textquotedblleft Pazy-Datko\textquotedblright\ result ( see Theorem 4.1, p. 116, of \cite{pazy}; also \cite{datko}, \cite{pazy2}), it is enough to establish the following integral estimate, for some $C^{\ast }>0$ that is independent of $t$: 
\begin{equation}
\int_{0}^{\infty }E(t)dt\leq C^{\ast }E(0).  \label{main}
\end{equation}%

At this point, we note that since $\Phi _{0}=\left[ p_{0},\mathbf{u}%
_{0},w_{0},w_{1}\right] \in \lbrack Null(\mathcal{A}_{0}\mathcal{)}]^{\bot },
$ then as pointed out in Remark \ref{inv}, $e^{\mathcal{A}_{0}t}\Phi _{0}\in
C([0,T];[Null(\mathcal{A)}]^{\bot })$ (See Proposition 13 of \cite%
{pelin-george}.) Now, with the objective estimate (\ref{main}) in mind, we
consider the elliptic map which was originally invoked in \cite{Chu2013-comp}
for the case $\mathbf{U=0.}$ Namely, let $\psi =\psi (f,g)\in H^{1}(\mathcal{%
O)}$ solve the following BVP for data $f\in L^{2}(\mathcal{O})$ and $g\in
L^{2}({\Omega }):$ 
\begin{equation}
\left\{ 
\begin{array}{c}
-\Delta \psi =f\text{ \ \ \ in \ }\mathcal{O} \\ 
\frac{\partial \psi }{\partial n}=0\text{ \ \ on \ }S \\ 
\frac{\partial \psi }{\partial n}=g\text{ \ \ on \ }\Omega 
\end{array}%
\right.   \label{156}
\end{equation}%
with $\left\{ f,g\right\} $ satisfying the compatibility condition%
\begin{equation*}
\int\limits_{\mathcal{O}}fd\mathcal{O+}\int\limits_{{\Omega }}gd{
\Omega=}0.
\end{equation*}%
By known elliptic regularity results for the Neumann problem on Lipschitz
domains--see e.g; \cite{JK}-- we have%
\begin{equation}
\left\Vert \psi (f,g)\right\Vert _{H^{\frac{3}{2}}(\mathcal{O)}}\leq C\left[
\left\Vert f\right\Vert _{\mathcal{O}}+\left\Vert g\right\Vert _{\Omega 
}\right] .  \label{156.5}
\end{equation}%
In order to get the estimate (\ref{main}), we invoke a multiplier born of $%
\psi (\cdot ,\cdot )$:\newline
\newline
\boxed{\textbf{Step I}:} With respect to the fluid PDE component in (\ref%
{k1}), we multiply both sides of this equation by $\nabla \psi (p(t),w(t)),$
where $\psi (p(t),w(t))$ satisfies (\ref{156}) with $f=p$ and $g=w$. Here,
we note that since $\Phi _{0}=\left[ p_{0},\mathbf{u}_{0},w_{0},w_{1}\right]
\in \lbrack Null(\mathcal{A)}]^{\bot }$, then $\{p(t),w(t)\}$ also satisfy
said compatibility condition for BVP (\ref{156}). An integration in time and
space then gives%
\begin{equation}
\int\limits_{0}^{T}\left( \mathbf{u}_{t}+\nabla p+\mathbf{U}\cdot \nabla 
\mathbf{u}-div~\sigma (\mathbf{u})+\eta \mathbf{u,}\nabla \psi (p,w)\right)
_{\mathcal{O}}=0.  \label{157}
\end{equation}%
We now\ look at each term in this relation separately; \newline
\boxed{\textbf{I.(i):}} We start with the first term of last relation and we have%
\begin{equation}
\int\limits_{0}^{T}(\mathbf{u}_{t}\mathbf{,}\nabla \psi (p,w))_{\mathcal{O}%
}d\tau =(\mathbf{u,}\nabla \psi (p,w))_{\mathcal{O}}|_{0}^{T}-\int%
\limits_{0}^{T}(\mathbf{u,}\nabla \psi (p_{t},w_{t}))_{\mathcal{O}}d\tau .
\label{158}
\end{equation}%
Now, let us focus on the second term of RHS of (\ref{158}): Using the
pressure equation in (\ref{k1})-(\ref{k3}), we get%
\begin{equation*}
-\int\limits_{0}^{T}(\mathbf{u,}\nabla \psi (p_{t},w_{t}))_{\mathcal{O}%
}d\tau =\int\limits_{0}^{T}(\mathbf{u,}\nabla \psi (\mathbf{U}\cdot \nabla
p,0))_{\mathcal{O}}d\tau 
\end{equation*}%
\begin{equation}
+\int\limits_{0}^{T}(\mathbf{u,}\nabla \psi (div~(\mathbf{u}),w_{t}))_{%
\mathcal{O}}d\tau .  \label{159}
\end{equation}%
Now, for the first term on RHS of (\ref{159}), we invoke the Leray (or
Helmholtz) Projector%
\begin{equation*}
\mathbb{P}\in \mathcal{L}\left( L^{2}(\mathcal{O)},L^{2}(\mathcal{O)}\cap
Null({div})\right) ;
\end{equation*}%
subsequently,%
\begin{equation*}
\mathbf{u}=\mathbb{P}\mathbf{u}+(I-\mathbb{P})\mathbf{u}
\end{equation*}%
satisfies%
\begin{equation}
\left\{ 
\begin{array}{c}
{div}(\mathbb{P}\mathbf{u})=0\text{ \ \ in \ }\mathcal{O};\text{ \ \ }%
\mathbb{P}\mathbf{u}\cdot \mathbf{n}=0\text{ \ \ on \ }\partial \mathcal{O}
\\ 
(I-\mathbb{P})\mathbf{u}=\nabla q(\mathbf{u}),\text{ \ \ \ for \ }q\in H^{1}(%
\mathcal{O)},\text{ \ }\int\limits_{\mathcal{O}}qd\mathcal{O=}0%
\end{array}%
\right.   \label{159.5}
\end{equation}%
(See e.g.; [Theorem 1.4, p. 11, \cite{temam}].) Therewith, we have%
\begin{equation*}
\int\limits_{0}^{T}(\mathbf{u,}\nabla \psi (\mathbf{U}\cdot \nabla p,0))_{%
\mathcal{O}}d\tau =\int\limits_{0}^{T}(\mathbb{P}\mathbf{u+}\nabla q(\mathbf{%
u})\mathbf{,}\nabla \psi (\mathbf{U}\cdot \nabla p,0))_{\mathcal{O}}d\tau 
\end{equation*}%
\begin{equation*}
=-\int\limits_{0}^{T}({div}\mathbb{P}\mathbf{u,}\psi (\mathbf{U}\cdot \nabla
p,0))_{\mathcal{O}}d\tau +\int\limits_{0}^{T}\left\langle \mathbb{P}\mathbf{%
u\cdot n,}\psi (\mathbf{U}\cdot \nabla p,0)\right\rangle _{\partial \mathcal{%
O}}d\tau 
\end{equation*}%
\begin{equation*}
+\int\limits_{0}^{T}\left\langle q(\mathbf{u})\mathbf{,}\nabla \psi (\mathbf{%
U}\cdot \nabla p,0)\cdot \mathbf{n}\right\rangle _{\partial \mathcal{O}%
}d\tau -\int\limits_{0}^{T}(q(\mathbf{u}){,{div}}\nabla \psi (\mathbf{U}%
\cdot \nabla p,0))_{\mathcal{O}}d\tau ;
\end{equation*}%
and so after considering (\ref{156}) and (\ref{159.5}) we have%
\begin{equation}
\int\limits_{0}^{T}(\mathbf{u,}\nabla \psi (\mathbf{U}\cdot \nabla p,0))_{%
\mathcal{O}}d\tau =\int\limits_{0}^{T}(q(\mathbf{u})\mathbf{,U}\cdot \nabla
p)_{\mathcal{O}}d\tau .  \label{160}
\end{equation}%
Now, to deal with RHS of (\ref{160}) we apply Green's Formula:%
\begin{equation*}
\int\limits_{0}^{T}(q(\mathbf{u})\mathbf{,U}\cdot \nabla p)_{\mathcal{O}%
}d\tau =\int\limits_{0}^{T}\int\limits_{\partial \mathcal{O}}(\mathbf{U}%
\cdot \mathbf{n})(q(\mathbf{u})p)d\partial \mathcal{O}
d\tau \end{equation*}%
\begin{equation*}
-\int\limits_{0}^{T}({div}(\mathbf{U)}q(\mathbf{u}),p)_{\mathcal{O}}d\tau
-\int\limits_{0}^{T}(\mathbf{U\cdot }\nabla q(\mathbf{u}),p)_{\mathcal{O}%
}d\tau 
\end{equation*}%
\begin{equation}
=-\int\limits_{0}^{T}({div}(\mathbf{U)}q(\mathbf{u}),p)_{\mathcal{O}}d\tau
-\int\limits_{0}^{T}(\mathbf{U\cdot (}[I-\mathbb{P}]\mathbf{u}),p)_{\mathcal{%
O}}d\tau .  \label{161}
\end{equation}%
Applying (\ref{159}) and (\ref{161}) to RHS of (\ref{158}) now gives 
\begin{equation*}
\int\limits_{0}^{T}(\mathbf{u}_{t}\mathbf{,}\nabla \psi (p,w))_{\mathcal{O}%
}d\tau =(\mathbf{u,}\nabla \psi (p,w))_{\mathcal{O}}|_{0}^{T}
\end{equation*}%
\begin{equation*}
-\int\limits_{0}^{T}({div}(\mathbf{U)}q(\mathbf{u}),p)_{\mathcal{O}}d\tau
-\int\limits_{0}^{T}(\mathbf{U\cdot (}[I-\mathbb{P}]\mathbf{u}),p)_{\mathcal{%
O}}d\tau 
\end{equation*}%
\begin{equation}
+\int\limits_{0}^{T}(\mathbf{u,}\nabla \psi (div~(\mathbf{u}),w_{t}))_{%
\mathcal{O}}d\tau .  \label{162}
\end{equation}%
In part, by using the classic estimate%
\begin{equation*}
\left\Vert q(u)\right\Vert _{\mathcal{O}}\leq C_{1}\left\Vert [I-\mathbb{P}%
]u\right\Vert _{\mathcal{O}}\leq C_{2}\left\Vert u\right\Vert _{\mathcal{O}},
\end{equation*}%
(see; e.g. Proposition 1.2, pp. 10 of \cite{temam} or Lemma 2.1.1(b), p.68
of \cite{sohr}), as well as the estimate for $\psi (\cdot ,\cdot )$ in (\ref%
{156.5})), we will then have from (\ref{162})%
\begin{equation}
\left\vert \int\limits_{0}^{T}(\mathbf{u}_{t}\mathbf{,}\nabla \psi (p,w))_{%
\mathcal{O}}d\tau \right\vert =\mathcal{O}\left( \Psi \left( \mathbf{U}%
\right) \left[ E(T)+E(0)\right] +\int\limits_{0}^{T}\left[ (\sigma (\mathbf{%
u(}\tau \mathbf{)}),\epsilon (\mathbf{u(}\tau \mathbf{))}_{\mathcal{O}}+\eta
\left\Vert \mathbf{u(}\tau \mathbf{)}\right\Vert _{\mathcal{O}}^{2}\right] ^{%
\frac{1}{2}}E(\tau )^{\frac{1}{2}}d\tau \right) ,  \label{163}
\end{equation}%
where we also implicitly use Korn's Inequality. (Here again, $\Psi \left( 
\mathbf{U}\right) \equiv \left[ \left\Vert \mathbf{U}\right\Vert _{\mathbf{H}%
^{3}(\mathcal{O})}+1\right] $). Now let us continue with the second term on
LHS of (\ref{157}).\newline
\boxed{\textbf{I.(ii):}} We have via Green's Formula and using (\ref{156}),%
\begin{equation*}
\int\limits_{0}^{T}(\nabla p\mathbf{,}\nabla \psi (p,w))_{\mathcal{O}}d\tau
=\int\limits_{0}^{T}\left\langle p\mathbf{n,}\nabla \psi (p,w)\right\rangle
_{\partial \mathcal{O}}d\tau -\int\limits_{0}^{T}(p\mathbf{,}\Delta \psi
(p,w))_{\mathcal{O}}d\tau 
\end{equation*}%
\begin{equation}
=\int\limits_{0}^{T}\left\Vert p\right\Vert _{\mathcal{O}}^{2}d\tau
+\int\limits_{0}^{T}\left\langle p\mathbf{n,}\nabla \psi (p,w)\right\rangle
_{\partial \mathcal{O}}d\tau .  \label{164}
\end{equation}%
\boxed{\textbf{I.(iii): }}To proceed with the fourth term on LHS of (\ref{157}), we
will appeal to the elliptic regularity results for solutions of second order
BVP on corner domains, which are established in \cite{Dauge_1,Dauge_2} (see
also \cite{Dauge_3}); it is at this point where our geometrical assumptions in Condition \ref{cond} come into play. Using these assumptions, we have the
following higher regularity,%
\begin{equation*}
\left\Vert \psi (p,w)\right\Vert _{H^{2}(\mathcal{O)}}\leq C\left[
\left\Vert p\right\Vert _{\mathcal{O}}+\left\Vert w_{ext}\right\Vert _{H^{%
\frac{1}{2}+\varepsilon }(\partial \mathcal{O)}}\right] 
\end{equation*}%
\begin{equation}
\leq C[\left\Vert p\right\Vert _{\mathcal{O}}+\left\Vert w\right\Vert
_{H_{0}^{2}(\Omega \mathcal{)}}],  \label{164.5}
\end{equation}%
where 
\begin{equation*}
w_{ext}(x)=\left\{ 
\begin{array}{c}
0,\text{ \ \ }x\in S \\ 
w(x),\text{ \ \ }x\in \Omega .%
\end{array}%
\right. 
\end{equation*}%
(For the second inequality in (\ref{164.5}), we are invoking Theorem 3.33,
pp. 95 of \cite{Mc}.) Therewith, for the fourth term on LHS of (\ref{157})
we have%
\begin{equation*}
-\int\limits_{0}^{T}(div~\sigma (\mathbf{u})\mathbf{,}\nabla \psi (p,w))_{%
\mathcal{O}}d\tau =\int\limits_{0}^{T}(\sigma (\mathbf{u})\mathbf{,}\epsilon
(\nabla \psi (p,w)))_{\mathcal{O}}d\tau -\int\limits_{0}^{T}\left\langle
\sigma (\mathbf{u})\cdot \mathbf{n,}\nabla \psi (p,w)_{\partial \mathcal{O}%
}\right\rangle d\tau 
\end{equation*}%
And so applying (\ref{164.5}) then gives%
\begin{equation*}
-\int\limits_{0}^{T}(div~\sigma (\mathbf{u})\mathbf{,}\nabla \psi (p,w))_{%
\mathcal{O}}d\tau 
\end{equation*}%
\begin{equation}
=\mathcal{O}\left( \int\limits_{0}^{T}\left[ (\sigma (\mathbf{u(}\tau 
\mathbf{)}),\epsilon (\mathbf{u(}\tau \mathbf{))}_{\mathcal{O}}+\eta
\left\Vert \mathbf{u(}\tau \mathbf{)}\right\Vert _{\mathcal{O}}^{2}\right] ^{%
\frac{1}{2}}E(\tau )^{\frac{1}{2}}d\tau \right)
-\int\limits_{0}^{T}\left\langle \sigma (\mathbf{u})\cdot \mathbf{n,}\nabla
\psi (p,w)_{\partial \mathcal{O}}\right\rangle d\tau   \label{165}
\end{equation}%
\newline
\boxed{\textbf{I.(iv):}} Lastly, by means of the regularity result in (\ref{156.5}%
), we have for the remaining terms in (\ref{157}), 
\begin{equation*}
\int\limits_{0}^{T}(\mathbf{U}\cdot \nabla \mathbf{u}+\eta \mathbf{u,}\nabla
\psi (p,w))_{\mathcal{O}}
\end{equation*}%
\begin{equation}
=\mathcal{O}\left( \int\limits_{0}^{T}\left[ (\sigma (\mathbf{u(}\tau 
\mathbf{)}),\epsilon (\mathbf{u(}\tau \mathbf{))}_{\mathcal{O}}+\eta
\left\Vert \mathbf{u(}\tau \mathbf{)}\right\Vert _{\mathcal{O}}^{2}\right] ^{%
\frac{1}{2}}E(\tau )^{\frac{1}{2}}d\tau \right) ,  \label{166}
\end{equation}%
where again we implicitly use Korn's ineqaulity. Now, if we combine (\ref%
{157}), (\ref{163}), (\ref{164}), (\ref{165}) and (\ref{166}), and keep in
mind that $[\sigma (\mathbf{u})\mathbf{n-}p\mathbf{n}]_{\partial \mathcal{O}%
}\cdot \mathbf{\tau }=0,$ as well as the BC in (\ref{156}) then we obtain%
\begin{equation}
\begin{array}{l}
\int\limits_{0}^{T}\left\Vert p\right\Vert _{\mathcal{O}}^{2}d\tau
-\int\limits_{0}^{T}\left\langle \left[ 2\nu \partial _{x_{3}}(\mathbf{u}%
)_{3}+\lambda \text{div}(\mathbf{u})-p\right] _{\Omega }\mathbf{,}%
w\right\rangle _{\Omega }d\tau  \\ 
\text{ \ \ \ \ \ \ \ \ \ \ \ \ }=\mathcal{O}\left( \Psi \left( \mathbf{U}%
\right) \left[ E(T)+E(0)\right] +\int\limits_{0}^{T}\left[ (\sigma (\mathbf{%
u(}\tau \mathbf{)}),\epsilon (\mathbf{u(}\tau \mathbf{))}_{\mathcal{O}}+\eta
\left\Vert \mathbf{u(}\tau \mathbf{)}\right\Vert _{\mathcal{O}}^{2}\right] ^{%
\frac{1}{2}}E(\tau )^{\frac{1}{2}}d\tau \right) .%
\end{array}
\label{167}
\end{equation}%
\boxed{\textbf{Step II:}} To continue with the energy estimates, in this
step we apply the multiplier $w$ to the plate equation in (\ref{k2}),
integrate in time and space to have%
\begin{equation*}
\int\limits_{0}^{T}(w_{tt}+\Delta ^{2}w+\left[ 2\nu \partial _{x_{3}}(%
\mathbf{u})_{3}+\lambda \text{div}(\mathbf{u})-p\right] _{\Omega
},w)_{\Omega }d\tau =0.
\end{equation*}%
An integration by parts then gives%
\begin{equation*}
\int\limits_{0}^{T}\left\Vert \Delta w\right\Vert _{\Omega }^{2}d\tau
+\int\limits_{0}^{T}(\left[ 2\nu \partial _{x_{3}}(\mathbf{u})_{3}+\lambda 
\text{div}(\mathbf{u})-p\right] _{\Omega },w)_{\Omega }d\tau 
\end{equation*}%
\begin{equation}
=\mathcal{O}\left( \int\limits_{0}^{T}\left\Vert w_{t}\right\Vert _{\Omega
}^{2}d\tau +E(T)+E(0)\right)   \label{168}
\end{equation}%
Using the domain criterion (A.v) (for $\kappa =0$) and the Sobolev
Trace Theorem, we also have%
\begin{equation}
\int\limits_{0}^{T}\left\Vert w_{t}\right\Vert _{\Omega }^{2}d\tau
=\int\limits_{0}^{T}\left\Vert u_{3}\right\Vert _{\Omega }^{2}d\tau =%
\mathcal{O}\left( \int\limits_{0}^{T}\left[ (\sigma (\mathbf{u(}\tau \mathbf{%
)}),\epsilon (\mathbf{u(}\tau \mathbf{))}_{\mathcal{O}}+\eta \left\Vert 
\mathbf{u(}\tau \mathbf{)}\right\Vert _{\mathcal{O}}^{2}\right] d\tau
\right) .  \label{169}
\end{equation}%
Now, adding the relations (\ref{167}), (\ref{168}) and (\ref{169}) we have
now%
\begin{eqnarray*}
\int\limits_{0}^{T}E(\tau )d\tau  &\leq &C\left( \Psi \left( \mathbf{U}%
\right) \left[ E(T)+E(0)\right] +\int\limits_{0}^{T}\left[ (\sigma (\mathbf{%
u(}\tau \mathbf{)}),\epsilon (\mathbf{u(}\tau \mathbf{))}_{\mathcal{O}}+\eta
\left\Vert \mathbf{u(}\tau \mathbf{)}\right\Vert _{\mathcal{O}}^{2}\right] ^{%
\frac{1}{2}}E(\tau )^{\frac{1}{2}}d\tau \right.  \\
&&\text{ \ \ \ \ \ \ \ \ \ \ \ }\left. +\int\limits_{0}^{T}\left[ (\sigma (%
\mathbf{u(}\tau \mathbf{)}),\epsilon (\mathbf{u(}\tau \mathbf{))}_{\mathcal{O%
}}+\eta \left\Vert \mathbf{u(}\tau \mathbf{)}\right\Vert _{\mathcal{O}}^{2}%
\right] d\tau \right) ,
\end{eqnarray*}%
whence we obtain after using Young Inequality%
\begin{equation*}
(1-\varepsilon )\int\limits_{0}^{T}E(\tau )d\tau \leq C_{0}\Psi \left( 
\mathbf{U}\right) \left[ E(T)+E(0)\right] +C_{\epsilon }\int\limits_{0}^{T}%
\left[ (\sigma (\mathbf{u(}\tau \mathbf{)}),\epsilon (\mathbf{u(}\tau 
\mathbf{))}_{\mathcal{O}}+\eta \left\Vert \mathbf{u(}\tau \mathbf{)}%
\right\Vert _{\mathcal{O}}^{2}\right] d\tau 
\end{equation*}%
where $C$ is a constant independent of $T>0$. To conclude the proof of
Theorem \ref{uniform}, we invoke the energy relation (%
\ref{E2}) and get%
\begin{equation}
(1-\varepsilon )\int\limits_{0}^{T}E(\tau )d\tau \leq 2CE(T)+C_{\varepsilon
}\int\limits_{0}^{T}\left[ (\sigma (\mathbf{u(}\tau \mathbf{)}),\epsilon (%
\mathbf{u(}\tau \mathbf{))}_{\mathcal{O}}+\eta \left\Vert \mathbf{u(}\tau 
\mathbf{)}\right\Vert _{\mathcal{O}}^{2}\right] d\tau .  \label{inter}
\end{equation}
Using the energy relation (\ref{E2}) once more time, we have%
\begin{equation*}
E(T)+\int\limits_{0}^{T}\left[ (\sigma (\mathbf{u(}\tau \mathbf{)}),\epsilon
(\mathbf{u(}\tau \mathbf{))}_{\mathcal{O}}+\eta \left\Vert \mathbf{u(}\tau 
\mathbf{)}\right\Vert _{\mathcal{O}}^{2}\right] d\tau \leq E(0)+\frac{1}{2}%
\left\vert \int\limits_{0}^{T}\int\limits_{\mathcal{O}}{div}(\mathbf{U})%
\left[ \left\vert p(\mathbf{\tau })\right\vert ^{2}+\left\vert \mathbf{%
u(\tau )}\right\vert ^{2}\right] d\tau d\mathcal{O}\right\vert .
\end{equation*}%
Applying this to (\ref{inter}), we then obtain%
\begin{equation*}
(1-\varepsilon )\int\limits_{0}^{T}E(\tau )d\tau \leq C_{1}\left( \Psi
\left( \mathbf{U}\right) \right) E(0)+C_{2}\Psi \left( \mathbf{U}\right)
\left\Vert {div}(\mathbf{U})\right\Vert _{\infty
}\int\limits_{0}^{T}E(\tau )d\tau .
\end{equation*}%
This gives the estimate (\ref{main}), for 
\begin{equation*}
\left\Vert {div}(\mathbf{U})\right\Vert _{\infty }<\frac{1}{2C_{2}\Psi
\left( \mathbf{U}\right) }.
\end{equation*}
and with subsequent $C^{\ast }=\frac{C_{1}\left( \Psi \left( \mathbf{U}%
\right) \right) }{2-\epsilon }$. This completes the proof of Theorem \ref%
{uniform}.

\section{Acknowledgement}

The author would like to thank the National Science Foundation, and
acknowledge her partial funding from NSF Grant DMS-1616425 and NSF Grant
DMS-1907823.

\end{document}